\RequirePackage{fix-cm}
\documentclass{svjour3}                     
\smartqed  
\RequirePackage{graphicx}
\RequirePackage{amsfonts}
\RequirePackage{amsmath}
\RequirePackage{amssymb}
\RequirePackage{tipx}
\RequirePackage{wrapfig}
\RequirePackage{ragged2e}
\RequirePackage{bm}

\providecommand{\e}[1]{\ensuremath{\times 10^{#1}}}

\begin{document}

\title{The moment-generating function of the log-normal distribution, how zero-entropy principle unveils an asymmetry under the reciprocal of an action.}


\author{
Yuri Heymann\        
}


\institute{Yuri Heymann \at
              \emph{Address in der Schweiz:} 3 rue Chandieu, 1202 Geneva, Switzerland \\
              \email{y.heymann@yahoo.com}             
 }

\date{Received: date / Accepted: date}

\maketitle

\begin{abstract}
The present manuscript is about application of It{\^o}'s calculus to the moment-generating function of the lognormal distribution. While Taylor expansion fails when applied to the moments of the lognormal due to divergence, various methods based on saddle-point approximation conjointly employed with integration methods have been proposed. By the Jensen's inequality, the MGF of the lognormal involves some convexity adjustment, which is one of the aspects under consideration thereof. A method based on zero-entropy principle is proposed part of this study, which deviations from the benchmark by infinitesimal epsilons is attributed to an asymmetry of the reciprocal.  As applied to systems carrying vibrating variables, the partial offset by the reciprocal of an action, is a principle meant to explain a variety of phenomena in fields such as quantum physics. 

\end{abstract}

\section{Introduction}

The moment-generating function (MGF) of a random variable $X$ is commonly expressed as $M(\theta) = \mathbb{E}({e}^{\theta X})$ by definition. By Taylor expansion of ${e}^{\theta X}$ centered on zero, we get $ \mathbb{E}({e}^{\theta X})=\mathbb{E}\left( 1 + \frac{\theta X}{1!} + \frac{\theta^2 X^2}{2!} + ... +  \frac{\theta^n X^n}{n!} \right)$. Taylor expansion method fails when applied to the MGF of the lognormal distribution of finite multiplicity. Its expression $M(\theta)= \sum_{n=0}^{\infty} \frac{\theta^n}{n!} {e}^{n \mu + n^2 \sigma^2 / 2}$ diverges for all $\theta$ - mainly because the lognormal distribution is skewed to the right. This skew increases the likelihood of occurrences departing from the central point of the Taylor expansion. These occurrences produce moments of higher order, which are not offset by the factorial of $n$ in the denominator, resulting in Taylor series to diverge. Some of the specificities of the lognormal distribution are enumerated below, namely that the lognormal distribution is not uniquely determined by its moments as seen in \cite{Heyde1963} for some multiplicity. Though it was reported in \cite{Romano} that MGF of the lognormal is not finite when $\theta$ is positive due to its integral form undefined for such $\theta$ values; other methods under consideration yield consistent values in a portion of the domain having positive $\theta$ values, e.g. Monte Carlo simulation, etc. 

\vspace{2mm}
\noindent
Saddle points and their representations are worth considering, though not as accurate in all domains. Expressions for the characteristic function of the lognormal distribution based on some classical saddle-point representations are described in \cite{Bruijn,Holgate1989}, and other representations involving asymptotic series in \cite{Barakat1976,Leipnik1991}. The evaluation of the characteristic function by integration method such as the Hermite-Gauss quadrature as in \cite{Gubner2006} appears to be the preferred choice. It was reported in \cite{beaulieu2004}, that due to oscillatory integrand and slow decay rate at the right tail of the lognormal density function, numerical integration was combersome. Other numerical methods including the work of \cite{Tellambura2010} employing a contour integral passing through the saddle point at steepest descent was proposed to overcome the aforementioned issues. 
A closed-form approximation of the Laplace transform of the lognormal distribution pretty accurate over its domain of definition is provided by \cite{Asmussen2016}. This equation is handy to backtest present study results, by its simplicity and because the Laplace transform of the density function and MGF are interconnected by sign interchange of variable $\theta$ (in the argument). An accurate and efficient method to value the lognormal MGF, referred to as thin-tile integration is proposed in section 4.0. This method which benefits from the symmetry of the Gaussian distribution is built on a partitioning of the density function into tiny elements having for basis a non-uniform grid spacing. 

\vspace{2mm}
\noindent
The approach consisting at using a stochastic process as a proxy of the MGF, is built upon the application of It{\^o}'s lemma with function $f(x) = {e}^{\theta Y(x)}$ to base process $dx_t = \mu dt + \sigma dW_t$, where $\mu$ and $\sigma$ are respectively the timeless mean and volatility parameters, $W_t$ a Wiener process, and $Y(x)$ is the random variable of the distribution resulting from application of some function to Gaussian variable $x$. As such the MGF is valued numerically by computation of $\mathbb{E}(f_t)=\mathbb{E}(f(x_t))$. In context of the MGF of the lognormal, the stochastic approach leads to tiny deviations explained by the partial offset  of the reciprocal of an action as the result of an asymmetry as applied to systems carrying vibrating variables. This is a principle having potential applications in other fields such as in quantum theory, e.g. asymmetries of harmonic oscillators \cite{Bruno1988,Jihad2020}. For example, the Heisenberg uncertainty principle implies the energy of a system described by harmonic oscillators cannot have zero energy \cite{Sciama1991}. The vibrational energy of ground state referred to as zero-point energy, is commonly invoked as a founding principle preventing liquid helium from freezing at atmospheric pressure regardless of temperature. 

\section{Theoretical background}

\subsection{The epsilon probability measure as a transform of zero convexity}

\noindent
By formal definition, a probability space is a measurable space satisfying measure properties such as countable additivity, further equipped with some probability measure assigning values to events of the probability space, e.g. value $0$ is assigned to the empty set and $1$ to the entire space. In field of stochastic calculus as applied to actuarial sciences, a probability measure is a mean to characterise a process by applying a drift to a base process resulting in an equivalence under the new measure satisfying some desirable properties e.g. is a martingale, etc. As such the epsilon probability measure belonging to probability space $\mathcal{E}_{\cdot}$ is equipped with a linear operator $\mathbb{E}_{\ast}$ as defined further down.

\vspace{2mm}
\noindent
Say $f\colon {\Bbb R} \to {\Bbb R}$ is a continuous function and $X$ a random variable, where $F$ is the primitive of $f$. Suppose there exists some measure in relation to $F(X)$ such that:

\begin{equation}
\mathbb{E}_{\ast} ( F(X) )= F(\mathbb{E}(X)) \,,
\end{equation}

\noindent
where $\mathbb{E}_{\ast}$ is an operator for statistical expectation of a variable in a probability space $\mathcal{E}_{Q}$ referring to such Q-measure, and where $\mathbb{E}$ is the corresponding operator in natural probability space. 

\vspace{2mm}
\noindent
As statistical expectation of a random variable can be expressed by its integral times density function over codomain in $\mathbb{R}$, the probability measure as defined in (1) expresses a kind of homeomorphism of the initial random variable by some dense function. Without specifications of higher moments, the primitive of $f$ as a diffeomorphism of multiplicity one i.e. stemming from a singleton implies there is an infinity of such variations spanning an entire domain or the probability space itself is empty. As such the $\mathcal{E}_{Q}$ probability space as a spectrum carrying operator $\mathbb{E}_{\ast}$ is generated by a one degree of freedom univariate, commonly characterised by mean of a class of functions or say a holomorphism expressing a multitude of moments by some multiplicity. 

\vspace{2mm}

\noindent
For example, if $F(X) = X^2$ where $X$ is a normally distributed variable centered in zero, then the above probability measure is not defined. As such for probability space $\mathcal{E}_{Q}$ to be non empty, implies $F$ is either monotonically increasing or decreasing over its domain, i.e. $F$ is a bijective application. 

\vspace{2mm}

\noindent
As a corollary stemming from (1), we can write:

\begin{equation}
\mathbb{E}_{\ast} \left( \int_{x_0}^{X} f(x) dx \right)= \int_{x_0}^{\mathbb{E}(X)} f(x) dx \,,
\end{equation}

\noindent
where $x_0$ is a real number and $X$ a random variable, yielding a new class of integrals represented as univariate functions, which integration domain is random on the right and braced to the left.

\vspace{1mm}

\noindent
A straight-forward property of statistical expectation as an operator belonging to probability space $\mathcal{E}_{Q}$ is linearity, yielding:

\begin{equation}
\mathbb{E}_{\ast}(a+ b X)= a + b  \, \mathbb{E}_{\ast}(X) \,,
\end{equation}

\noindent
where $a$ and $b$ are real numbers and $X$ a variable carrying some randomness. This relation is viewed as a diffeomorphism of the density function by mean of variations spanning a quantised space, stemming from integral representation of statistical expectation by a dense function of $X$ in new probability space $\mathcal{E}_{Q}$, where $Q$ stands as measure. 

\vspace{2mm}

%
%

\noindent{\textbf{Proposition 1 - applies to Gaussian scenerio:}} Given an epsilon probability measure in relation to function $f$ which applies to time-sensitive variable $X_t$ linked to a Wiener process on a one-by-one relation where $f(X_t)$ is Gaussian, we have:

\begin{equation}
\mathbb{E}_{\ast}(f(X_t)) =  \mathbb{E}(f(X_t))\,,
\end{equation}

\noindent
where operator $\mathbb{E}_{\ast}$ belonging to epsilon probability space $\mathcal{E}_{Q}$ is linear and event space remains Hausdorff under the Q-measure. 

\begin{proof}
Say we have univariate function $f$ which applies to variable $X_t$ of a random nature such that $f(X_t)$ is Gaussian. For process $f(X_t)$ as a Gaussian variable sensitive to time $t$ and equipped with a Wiener process $W_t$, where there is a bijective map between $X_t$ and $W_t$, the latter implies the existence of two real numbers $a$ and $b$ such that $f(X_t)=a + b \,W_t$. As $\mathbb{E}_{\ast}(W_t)=0$ in $\mathcal{E}_{Q}$ (by linearity of the operator), we have $\mathbb{E}_{\ast}(f(X_t))= a = \mathbb{E}(f(X_t))$ which constitutes proof of the above.
\end{proof}

\subsection{Fundamentals related to an action and reciprocal for systems carrying vibrating variables}


\noindent
Say we apply a function $f \colon \mathbb{R} \to \mathbb{R}$ to a vibrating variable $x_t$, which action is projecting $x_t$ onto $f(x_t)$ in a measurable space expressing statistical expectations as a linear operator. The action resulting from the application of It{\^o}'s lemma to variable $x_t$ of a base process with function $f$, is leading to derived process $f_t=f(x_t)$. The reciprocal as a reverse transform is aiming at projecting the statistical expectations of $f(x_t)$ onto $x_t$ by linear operator, which in context of the lognormal distribution and some normalization leads to some skewness.

\vspace{2mm}
\noindent
By the Jensen's inequality, the application of a function $f \colon \mathbb{R} \to \mathbb{R}$ to a vibrating variable $x_t$ yields a convexity adjustment $\Delta f_t$, by relation $\mathbb{E}(f(x_t))= f(\mathbb{E}(x_t)) + \Delta f_t$ expressing a projection of statistical expectation under a skewed distribution. The shifted process $\widetilde{x_t}$ is an implied variable of $x_t$ which by value satisfies equality $\mathbb{E}(f(\widetilde{x_t}))= f(\mathbb{E}(x_t))$, as a link with epsilon probability measure introduced earlier.

\vspace{2mm}
\noindent
For such a primary transform projecting the statistical expectation of $x_t$ onto $f(\widetilde{x_t})$ by a shift linked to measurable space epsilon, the reciprocal obtained from the application of the inverse function $f^{-1}$ to process of $f(\widetilde{x_t})$  results in projection of $f(\widetilde{x_t})$ onto $x_t$ by a linear operator. 

\vspace{2mm}
\noindent
Given a base process $dx_t=\mu dt + \sigma dW_t$ where $\mu$ and $\sigma$ are the timeless mean and volatility parameters and $W_t$ a Wiener process, the It{\^o}'s process of $f_t$ is defined such that $f_t=f(x_t)$ where function $f\colon {\Bbb R} \to {\Bbb R}$ is a twice-differentiable function, continuous to the right and invertible. The It{\^o}'s process of $f_t$ obtained by applying It{\^o}'s lemma to base process $x_t$ with function $f$, leads to:

\begin{equation}
df_t = \left( \mu \frac{\partial f (x_t)}{\partial x} + \frac{\sigma^2}{2} \frac{\partial^2 f (x_t)}{\partial x^2} \right) dt + \sigma \frac{\partial f}{\partial x}(x_t) \, dW_t\,,
\end{equation}

\noindent 
where $W_t$ is a Wiener process, $\mu$ and $\sigma$ are real numbers representing the timeless mean and volatility of base process $x_t$, where $t$ is time.

\vspace{2mm}
\noindent
When normalizing the It{\^o} process of $f_t$, by dividing both sides of (5) by partial derivative $\partial f / \partial x$ of $x_t$, leads to stochastic differential equation (SDE):

\begin{equation}
[f^{-1}]^{\prime}(f_t) df_t= \left( \mu + \frac{1}{2}\sigma^2 \, \frac{\partial^2 f (x_t)}{\partial x^2} \middle/ \frac{\partial f (x_t)}{\partial x}  \right) dt + \sigma \, dW_t \,,
\end{equation}

\noindent
where $[f^{-1}]^{\prime}(f_t)$ on the left-hand side of (6) is coming from relation $[f^{-1}]^{\prime}(a)=1/f^{\prime}(f^{-1}(a))$ per the derivative of the inverse of $f$ evaluated in $a$ real.

\vspace{2mm}
\noindent
By applying It{\^o}'s lemma to the It{\^o}'s process of $f_t$ given by (5) with function $f^{-1}$ yields:

\begin{equation}
d(f^{-1}(f_t))= [f^{-1}]^{\prime}(f_t) \, df_t - h(f_t) \, dt \,,
\end{equation}

\noindent
where $h(f_t)$ is expressed as follows:

\begin{equation}
h(f_t) = h(f(x_t))= \frac{1}{2} \sigma^2 \frac{f^{\prime\prime}(x_t)}{f^{\prime}(x_t)} \,.
\end{equation}

\noindent
The term $h(f_t) dt$ in (7) is resulting from the convexity adjustment of transform $T \colon f_t \to f^{-1}(f_t)$ as a projection of statistical expectation by application of It{\^o}'s lemma to base process and normalisation as seen above. 
Under perfect symmetry, we would expect the below to happen:

\vspace{2mm}
\noindent
\noindent{\textbf{(i)}} By removing the vibrations in (7) conjointly with the application of linear operator $\mathbb{E}_{\ast}$ to the integral form of the stochastic process (6) as seen in (2), leads to a shift of $f_t$ towards $\widetilde{f_t}$ under the action of $T$. More specifically, $\widetilde{f_t}$ is the implied process of $f_t$ which by value satisfies equality $d(f^{-1}(\widetilde{f_t}))= [f^{-1}]^{\prime}(\widetilde{f_t}) \, d\widetilde{f_t}$ also referred to as the zero-entropy condition. The combined action of the above leads to transform $T_{\ast} \colon f_t \to f^{-1}(\widetilde{f_t})$ as a projection of statistical expectations.

\vspace{2mm}
\noindent
\noindent{\textbf{(ii)}} As transform $T_{\ast} \colon f_t \to f^{-1}(\widetilde{f_t})$ represents a projection of statistical expectations through a linear operator, by the symmetry of the reciprocal we have $f^{-1}\left(\mathbb{E}(f_t)\right) = \mathbb{E} \left(f^{-1}(\widetilde{f_t})\right)$. The benefit sought by such symmetry is a simple expression for $\mathbb{E}(f_t)$ which is equal to $f(\mathbb{E}(f^{-1}(\widetilde{f_t})))$ under perfect symmetry. Yet, this is not always as accurate due to partial offset of the shifts resulting from the action of $f^{-1}$ and reciprocal, which asymmetry is explained by non-offseting convexity adjustments. 

\vspace{2mm}
\noindent
From the linearity of statistical operators applicable to stochastic processes, we expect the reciprocal of an action to offset the action of the primary transform, which can be easily verified for MFG of Gaussian random variables, expressing a projection of statistical expectations from epsilon space to natural probabilities. 

\vspace{2mm}
\noindent
A variant of the above is to say that $x_t \sim \mathcal{N}(\mu \, t, \sigma^2 \,t)$ is a Gaussian random variable, where $\mu \, t$ is the mean and $\sigma^2 \, t$ the variance. Given a bijective function $f$ of real codomain and It{\^o} process $f_t$ defined as $f_t =f(x_t)$, where $f^{-1}(f_t)$ has Gaussian distribution $\mathcal{N}(\mu_t, v_t)$ with mean $\mu_t$ and variance $v_t$, the above is equivalent to finding the implied convexity of function $f^{-1}$ as applied to $f_t$ for any $t \geq 0$. The trivial case when $f$ is a linear function leads to complete offset of the effect of the vibrations of $x_t$ by the reciprocal. In this scenario, variances before applying transformation $f^{-1}$ and after reciprocal $f$ is applied are identical. Yet, the approach consisting at computing the implied convexity from the variance of the reciprocal of an arbitrary transform remains difficult in practice. A rather suitable method referred to as thin-tile integration, which is based on non-uniform grid spacings is provided further down, see §4.0. The proceeding is the continuation of the former stochastic approach to the lognormal distribution.

\subsection{The zero-entropy principle and Gaussianity of stochastic processes resulting from non-linear transformations of a normally distributed process}

\noindent
By considering a base process $x_t$  having distribution $\mathcal{N}(\mu_1 \, t, \sigma_1^2 \, t)$, its  SDE is expressed as:
 
\begin{equation}
dx_t = \mu_1 dt + \sigma_1 \, dW_t\,,
\end{equation}

\noindent
where $W_t$ is a Wiener process, $\mu_1$ and $\sigma_1$ real numbers representing timeless mean and volatility parameters.

\vspace{2mm}
\noindent
Say $f$ is a twice-differentiable function continuous to the right and invertible defined as follows:

\begin{equation}
f(x)= \int_{x_0}^{x} \sigma_2(s) ds \,,
\end{equation}

\noindent
where $\sigma_2$ is a univariate function in $\mathbb{R}$ and $x_0$ the lower bound of the integration domain.

\vspace{2mm}
\noindent
Given $\frac{\partial f(x)}{\partial x} = \sigma_2(x) $ and $\frac{\partial^2 f(x)}{\partial x^2} =  \sigma_2^{\prime}(x)$, the application of It{\^o}'s lemma to the base process (9) with function in (10), leads to:

\begin{equation}
df_t = \left( \mu_1 \sigma_2(f^{-1}(f_t)) + \frac{1}{2}\sigma_1^2 \sigma_2^{\prime}(f^{-1}(f_t)) \right)dt + \sigma_1 \sigma_2(f^{-1}(f_t)) \, dW_t \,,
\end{equation}

\noindent
where $f^{-1}$ is the inverse of $f$ and $f_t$ the It{\^o}'s process defined as $f_t=f(x_t)$, where $\sigma_2(f^{-1}(f_t))= \sigma_2(x_t)$. By the normalization of (11) with $\sigma_2(x_t)$, leads to:

\begin{equation}
\frac{df_t}{\sigma_2(x_t)} = \left( \mu_1 + \frac{1}{2}\sigma_1^2 \, \frac{ \sigma_2^{\prime}(x_t)}{\sigma_2(x_t)} \right) dt + \sigma_1 \, dW_t \,,
\end{equation}

 \noindent
By the aggregation of the drift in (12) into single expression $\mu_2(x_t)$, yields:

\begin{equation}
\int_{f_0}^{f_t} \frac{df}{\sigma_2( f^{-1}(f))} = \int_{0}^{t}\frac{\mu_2(f^{-1}(f_s))}{\sigma_2(f^{-1}(f_s))} ds + \int_{0}^{t} \sigma_1 dW_t \,,
\end{equation}

\noindent
which is the corresponding integral form of SDE (12). By the derivative of the inverse of a function, say in a zero-entropy universe, we have: 

\begin{equation}
 f^{-1}(\xi) = \int_{f_0}^{\xi} \frac{ds}{\sigma_2( f^{-1}(s))}\,,
\end{equation}

\noindent
where $\xi$ is a non-vibrating variable. The application of It{\^o}'s lemma to the process of $f_t$ with function $f^{-1}(\xi)$, would break equality $f^{-1}(f_t) = \int_{f_0}^{f_t} \frac{df_t}{\sigma_2( f^{-1}(f_t))}$ as defined in (14) for zero-entropy, i.e. where $\xi$ is a non-vibrating variable. Say $\widetilde{f_t}$ is the shifted process of $f_t$ satisfying the equality for zero-entropy to hold. Thus, we have:
\begin{equation}
f^{-1}(\widetilde{f_t})=  \int_{0}^{t}\frac{\mu_2(f^{-1}(\widetilde{f_s}))}{\sigma_2(f^{-1}(\widetilde{f_s}))} ds + \sigma_1 W_t \,,
\end{equation}

\newpage
\noindent{\textbf{Proposition 2 - Gaussianity condition:}} If $f^{-1}$ spans $\mathbb{R}$ on the image of $f$, then $f^{-1}(\widetilde{f_t})$ is Gaussian with distribution $\mathcal{N}(m_t, v_t)$ where $m_t$ is the mean and $v_t$ the variance of the process by filtration $\mathcal{F}_t$  at any instant of time $t \geq 0$. As a counterexample, if function $f(x)=\sqrt{x}$ is used in (10) where $x \in \mathbb{R}^+$, then $f^{-1}(\widetilde{f_t})$ departs from Gaussianity when approaching the left part of the domain. 

\begin{proof}
Consider a stochastic process $X_t$ defined such that $X_t = \int_{0}^{t}g(X_s) ds + \sigma W_t$ where $W_t$ is a Wiener process, $\sigma$ timeless volatility parameter and $g$ is a real-valued function continuous in  $\mathbb{R}$; thus: $dX_t=g(X_t) dt+\sigma dW_t$. As $X_t$ can be decomposed into an infinite sum of tiny independent Gaussian random variables, implies that $X_t$ is Gaussian. By considering a discretization of the process into infinitesimal time intervals $\delta t$, we have $X_{t+1}-X_t = g(X_t) \delta t +\sigma \sqrt{\delta t} Z_t$ where $Z_t$ is a standard normal random variable. By the first iteration,  $X_1 = X_0 +g(X_0) \delta t + \sigma \sqrt{\delta t} Z_0$. Corresponding increment expressed as $\delta X_0=g(X_0) \delta t + \sigma \sqrt{\delta t} Z_0$, is a Gaussian infinitesimal. By the second iteration, $X_2= X_1+g(X_0+ \delta X_0) \delta t + \sigma \sqrt{\delta t} Z_1$, leads to corresponding increment $\delta X_1=g(X_0+ \delta X_0) \delta t + \sigma \sqrt{\delta t} Z_1$. As increment $\delta X_0$ is infinitesimal, Taylor expansion of $g(X_0+ \delta X_0)$ leads to $g(X_0+ \delta X_0) \doteq g(X_0) + g^{\prime}(X_0) \delta X_0$, as higher order terms are negligible. Hence $\delta X_1=\left(g(X_0) + g^{\prime}(X_0) \delta X_0 \right) \delta t + \sigma \sqrt{\delta t} Z_1$. When $\delta t \rightarrow 0$, $\delta X_1$ tends to a Gaussian infinitesimal. By the third iteration: $X_3= X_2+g(X_0+ \delta X_0 + \delta X_1) \delta t + \sigma \sqrt{\delta t} Z_2$, with corresponding increment $\delta X_2=g(X_0+ \delta X_0+\delta X_1) \delta t + \sigma \sqrt{\delta t} Z_2$. Because $\delta X_0$ and $\delta X_1$ are infinitesimal increments, we expand $g(X_0+ \delta X_0+\delta X_1)$ by bivariate Taylor expansion applied to expression $\varphi(x, y)=g(X_0+ x+y)$ in $(x,y)=(\delta X_0, \delta X_1)$ to the first order (i.e. higher order terms negligible). We get $g(X_0+ \delta X_0 + \delta X_1) \doteq \varphi(0,0)+\varphi_x(0,0) \delta X_0 + \varphi_y(0,0) \delta X_1 = g(X_0) + g^{\prime}(X_0) \left(\delta X_0 + \delta X_1 \right)$. Hence, $\delta X_2=\left( g(X_0) +g^{\prime}(X_0) \left(\delta X_0 + \delta X_1 \right) \right) \delta t + \sigma \sqrt{\delta t} Z_2$. When $\delta t \rightarrow 0$, $\delta X_2$ tends to a Gaussian infinitesimal. By the successive application of the above steps, we get that for all integers $i=0,...,n$, $\delta X_i$ is a linear combination of $\delta X_0$,..,$\delta X_{i-1}$ and $Z_i$. By substitution of $\delta X_0$ into expression of $\delta X_1$ and so on, yieldsilt
expression of $\delta X_i$ as a linear combination of $Z_0$,..,$Z_i$ for all integers $i=0,...,n$. When adding together all the $\delta X_i$ sub-components, leads to an expression of $X_t$ as a linear combination of $Z_0$,...,$Z_n$, which are independent standard normal random variables. Hence, we can say that $X_t$ is Gaussian. In (15), $X_t=f^{-1}(\widetilde{f_t})$. A prerequisite for $X_t$ to be Gaussian, is that $f^{-1}$ spans the entire domain in $\mathbb{R}$ on the image of $f$ as given by $f_t=f(x_t)$. This constitutes proof of proposition 2. 
\end{proof}

\noindent
If $f^{-1}(\widetilde{f_t})$ is Gaussian (see \textit{proposition 2}), then linear operator $\mathbb{E}_{\ast}$ of the epsilon probability measure $\mathcal{E}_{Q}$ as defined earlier can be applied to $f^{-1}(\widetilde{f_t})$, which by \textit{proposition 1} leads to:

\begin{equation}
 \mathbb{E}_{\ast}(f^{-1}(\widetilde{f_t})) =  \mathbb{E}(f^{-1}(\widetilde{f_t})) = m_t \,.
\end{equation}

\noindent
Eq. (16) requires $\mathbb{E}_{\ast}(W_t) = 0$, which as a linear operator is true in a probability space, where the event space remains Hausdorff under the change of measure. This is a variant of sigma-fields, having an event space equipped of a linear operator, that can be partitioned into infinitely many collections of disjoint subelements which parts are subadditive, meaning outcomes are frictionless.  

\vspace{2mm}
\noindent
By perfect symmetry of an action and reciprocal under zero-entropy principle, we have $\mathbb{E}(f^{-1}(\widetilde{f_t})) =  f^{-1}(\mathbb{E}(f_t))$, leading to:

\begin{equation}
 f^{-1}(\mathbb{E}(f_t)) =  m_t \,,
\end{equation}

\noindent
as a prerequisite for the estimator of $\mathbb{E}(f(x_t))$ expressed as:

\begin{equation}
\hat{\mathbb{E}}(f_t)= f(m_t) \,,
\end{equation}

\noindent
where $\hat{\mathbb{E}}(f_t)$ is the corresponding estimator for statistical expectation of function $f$ as applied to variable $x_t$ of base process. 

\section{Characterisation of stochastic processes related to the lognormal MGF by zero-entropy principle} 

Let's apply the above to the MGF of the lognormal. By definition the MGF of the lognormal of parameter $\mu$ and $\sigma$ is expressed as $M(\theta) =\mathbb{E}({e}^{\theta {e}^x})$, where $x \sim \mathcal{N}(\mu, \sigma^2)$ and $\theta$ is real.

\subsection{From a-level stochastic calculus}

\noindent
Consider the base process $x_t$ of Gaussian distribution $x_t \sim \mathcal{N}(\mu t, \sigma^2 t)$, which is defined as:

\begin{equation}
dx_t = \mu dt + \sigma dW_t\,,
\end{equation}

\noindent
where $W_t$ is a Wiener process, $\mu$ and $\sigma$ the timeless mean and volatility parameters, and initial value given by $x_0=0$. 

\vspace{2mm}
\noindent
By applying It{\^o}'s lemma to (19) with $f(x)={e}^{\theta {e}^x}$ and normalisation as per (11-12), leads to:

\begin{equation}
\frac{df_t}{f_t \ln f_t} = \left( \mu +\frac{1}{2}\sigma^2 +\frac{1}{2}\sigma^2 \ln f_t \right) dt + \sigma dW_t\,.
\end{equation}

\noindent
From zero-entropy principle, we say there exists a shifted process $\widetilde{f_t}$ such that equality $d(\ln \ln \widetilde{f_t})= \frac{d\widetilde{f_t}}{\widetilde{f_t} \ln \widetilde{f_t}}$ holds, leading to:

\begin{equation}
\int_{f_0}^{\widetilde{f_t}} \frac{ds}{s \ln s} =\int_{0}^{t} \left( \mu +\frac{1}{2}\sigma^2 +\frac{1}{2}\sigma^2 \ln \widetilde{f_s} \right) ds + \sigma W_t\,.
\end{equation}

\noindent
By splitting its domain of definition into positive and negative parts, yields the below expressions resulting from sign interchange of $\theta$ through complex logarithm and reflection, i.e. when $\theta \in \mathbb{R}^{+}$:

\begin{equation}
\ln \ln (\widetilde{f_t})  =\ln \theta + \left( \mu +\frac{1}{2}\sigma^2 \right) t + \frac{1}{2}\sigma^2 \int_{0}^{t} \ln (\widetilde{f_s})ds + \sigma W_t \,,
\end{equation}

\noindent
and when $\theta \in \mathbb{R}^{-}$:

\begin{equation}
\ln (-\ln (\widetilde{f_t}))  =\ln (-\theta)+ \left( \mu +\frac{1}{2}\sigma^2 \right) t - \frac{1}{2}\sigma^2 \int_{0}^{t} \ln (\widetilde{f_s})ds + \sigma W_t \,.
\end{equation}

\vspace{2mm}
\noindent
As per \textit{proposition 2}, we can say that for all $\theta$ positive, (22) is a $\mathcal{F}_t$-adapted Gaussian process with distribution $\mathcal{N}(m_t, v_t)$ where $m_t$ is the first moment and $v_t$ the variance of the process at any instant of time $t \geq 0$. The same can be said about Gaussianity of $\ln (-\ln (\widetilde{f_t}))$ in (23) when $\theta$ negative. To handle both parts of the domain into a single expression, the sign function denoted $\text{sign}(x)$ (which returns $1$ when $x$ is positive and $-1$ when negative), is invoked further down in the remaining of the manuscript.

\subsection{The underpinning structure between the real and imaginary parts of process $y_t$ when extended to complex numbers}

\noindent
By setting $y_t=\ln (\text{sign}(\theta) \ln \widetilde{f_t})$, (22-23) can be rewritten as follows:

\begin{equation}
y_t  =\ln (|\theta|)+ \left( \mu +\frac{1}{2}\sigma^2 \right) t + \text{sign}(\theta) \, \frac{1}{2}\sigma^2 \int_{0}^{t} {e}^{y_s} ds + \sigma W_t \,.
\end{equation}

\noindent
As an extension to the complex domain, the above process $y_t$ can be split into real and imaginary components according to $y_t = y_{R,t}+i \, y_{I,t}$, where $y_{R,t}$ and $y_{I,t}$ are the real and imaginary parts respectively. In most settings, the bivariate normal distribution relies on law of large numbers, a particular case of multivariate analysis. As from (24), the dependence between the marginals of $y_t$ is as follows:

\begin{equation}
y_{R, t} = \Re(\ln |\theta|) + \left(\mu+\frac{1}{2} \sigma^2\right) t + \text{sign}(\theta) \,\frac{1}{2} \sigma^2 \int_{0}^{t} \cos (y_{I, s}) {e}^{y_{R, s}} ds + \sigma W_t \,,
\end{equation}

\noindent
and

 \begin{equation}
y_{I, t} = \Im (\ln |\theta|) +\text{sign}(\theta) \, \frac{1}{2} \sigma^2 \int_{0}^{t} \sin (y_{I, s}) {e}^{y_{R, s}} ds\,.
\end{equation}

\noindent
where the logarithm is extended to $\mathbb{C}$ by inversion of complex exponentiation, carrying signs as per the complex quadrant, and where $\text{sign}(\theta)$ has correspondance with sign of $\Re(\theta)$ (see de Moivre).

\vspace{2mm}
\noindent
Despite the real part of process $y_t$ as an $\mathcal{F}_t$-adapted process is Gaussian, its imaginary counterpart $y_{I, t}$ is not $\mathcal{F}_t$-adapted as it carries some kind of dependence on the real part of a path dependent nature. Though carrying some correlation with real part $y_{R, t}$ vibrating according to a Gaussian spread, imaginary part $y_{I,t}$ exhibits non-Gaussianity as the result of the application of Euler's formula to $e^{i \, y_{I,t}}$, leading to a 
diffusive function of $y_{R, t}$ by filtration $\mathcal{F}_t$ respective to Wiener process $W_t$, see (26). It follows that $y_{I, t}$ is non-Gaussian by non-linearity respective to $\mathcal{F}_{t}$-adapted Gaussian process $y_{R, t}$ in its probability measure, and as a composite translated into bidimensional-adapted process by say filtration $\mathcal{F}_{t,s}$, where $s$ is path dependent dimensionality acting as dispersive factor carrying some skew by some non-linear operator.

\subsection{Expression to characterise the first moment of process $y_t$ at time $t$}

\noindent
Rewriting (22-23) as:

\begin{equation}
y_t  =\ln (|\theta|) + \left( \mu +\frac{1}{2}\sigma^2 \right) t +\text{sign}(\theta) \, \frac{1}{2}\sigma^2 \int_{0}^{t} {e}^{y_s} ds + \sigma W_t \,,
\end{equation}

\noindent
where $y_t=\ln (\text{sign}(\theta) \,\ln \widetilde{f_t})$ and $\theta$ is real, implies by \textit{proposition 2} that as an $\mathcal{F}_t$-adapted process, $y_{t}$ is Gaussian with mean $m_t$ and variance $v_t$, i.e. $y_t \sim \mathcal{N}(m_t, v_t)$.

\vspace{2mm}
\noindent
By applying the statistical expectation as a linear operator to SDE (27) and using relation $\mathbb{E}(\int_{0}^{t} h(s) ds) = \int_{0}^{t} \mathbb{E} \left( h(s) \right) ds$ resulting from linearity, we get:

\begin{equation}
\mathbb{E} \left( y_t \right)  =\ln (|\theta|) + \left( \mu +\frac{1}{2}\sigma^2 \right) t + \text{sign}(\theta) \, \frac{1}{2}\sigma^2 \int_{0}^{t} \mathbb{E} \left( {e}^{y_t} \right) ds \,.
\end{equation}

\noindent
The expected value of a lognormal variable expressed as $z={e}^x$, where $x \sim \mathcal{N}(\mu, \sigma^2)$ equals $\mathbb{E}(z) = {e}^{\mu + \frac{1}{2}\sigma^2}$, thus leading to:

\begin{equation}
\mathbb{E} \left( {e}^{y_t} \right)  = {e}^{(m_t+ \frac{1}{2} v_t)} \,.
\end{equation}

\noindent
As a note (29) holds when $\theta$ is a real, say for the lognormal MGF when $\theta$ has negative or positive values see (22-23). This equality is no longer true when $\theta$ is extended to the complex domain, as only the real component of $y_t = y_{R,t}+i \, y_{I,t}$ is Gaussian. Due to non-Gaussianity of the imaginary part $y_{I, t}$ as seen in §3.2, process $y_t$ is not an $\mathcal{F}_t$-adapted Gaussian process. The linear combination of a Gaussian and non-Gaussian variable yielding a non-Gaussian process, implies that Gaussianity of $y_t$ cannot be invoked for the statistical expectation $\mathbb{E} \left( {e}^{y_t} \right)$.

\vspace{1mm}
\noindent
We can rewrite (28) as follows:

\begin{equation}
m_t  = \ln (|\theta| ) + \left( \mu +\frac{1}{2}\sigma^2 \right) t +\text{sign}(\theta) \, \frac{1}{2}\sigma^2 \int_{0}^{t} {e}^{(m_s+ \frac{1}{2} v_s)} ds \,.
\end{equation}

 \noindent
As $m_t$ and $v_t$ can be viewed as real-valued functions sensitive to time $t$, no rules prevents direct derivation of expression (30) by the $\partial / \partial t$ operator, leading to:

\begin{equation}
\frac{ \partial m_t}{\partial t}  = \mu +\frac{1}{2}\sigma^2  + \text{sign}(\theta) \,\frac{1}{2}\sigma^2  {e}^{(m_t+ \frac{1}{2} v_t)}  \,,
\end{equation}

\noindent
with initial condition $m_0= \ln(|\theta|)$. Eq. (31) is the differential equation describing the first moment of process $y_t$, at any instant of time $t \geq 0$ and where $\theta$ is real.

\subsection{Expression to characterise the variance of process $y_t$ at arbitrary time}

\noindent
As an $\mathcal{F}_t$-adapted process where $\theta$ is real, (27) is Gaussian i.e. $y_t \sim \mathcal{N}(m_t, v_t)$, which in affine representation can be rewritten as:

\begin{equation}
\widetilde{y_t} = m_t + \sqrt{v_t} Z\,,
\end{equation}

\noindent
where $\widetilde{y_t}$ is the affine representation of $y_t$, with $t \geq 0$ representing time and standard Gaussian by $Z \sim \mathcal{N}(0, 1)$ . 

\vspace{2mm}
\noindent
As process in (32) is expressed in pure notations, we can directly apply operator $\partial / \partial t$ as a partial derivative to component functions, leading to:

\begin{equation}
\frac{\partial \widetilde{y_t}}{\partial t} = m_t^{\prime} + \frac{1}{2} \frac{v_t^{\prime}}{\sqrt{v_t}} Z\,,
\end{equation}

\noindent
where $m_t^{\prime}$ and $v_t^{\prime}$ are the respective time derivatives of $m_t$ and $v_t$ as functions. The variance of (33) is further expressed as:

\begin{equation}
Var\left( \frac{\partial \widetilde{y_t}}{\partial t} \right) = \frac{1}{4} \frac{{v_t^{\prime}}^2}{v_t} \,.
\end{equation}

\noindent
As $y_t=\ln (\text{sign}(\theta) \, \ln \widetilde{f_t})$, SDE (27) expressed in infinitesimal form is as follows:

\begin{equation}
dy_t = \left( \mu +\frac{1}{2}\sigma^2 + \text{sign}(\theta) \, \frac{1}{2}\sigma^2 {e}^{y_t} \right) dt + \sigma dW_t\,,
\end{equation}

\noindent
As time integration viewed as an operator applied to an It{\^o}'s processes expressing a simple difference $dy_t$ on the left-hand side of the SDE is not convex sensitive by Gaussianity of $y_t$, integration of (35) with respect to time, leads to:

\begin{equation}
y_t - y_0 =  \int_{0}^{t} \left( \mu +\frac{1}{2}\sigma^2 +\text{sign}(\theta) \,\frac{1}{2}\sigma^2 {e}^{y_t}  \right) dt + \sigma W_t \,.
\end{equation}

\vspace{2mm}
\noindent
The application of operator $\partial / \partial t$ as a partial derivative to functions carrying vibrations as in (36) where the Wiener process is in its aggregated form $W_t=\sqrt{t} \, Z$ with $Z \sim \mathcal{N}(0, 1)$ requires special attention. Say SDE in (36) is expressed as $y_t=f(y_t) + x_t$, where $x_t=\sigma \sqrt{t} Z$ is the source of vibrations. 

\vspace{3mm}
\noindent
The relation between $y_t$ and $x_t$ is obained from bilinearity between affine Gaussian representation $y_t = a_t + b_t \, Z$ and $x_t = \sigma \sqrt{t} Z$, leading to:

\begin{equation}
y_t = a_t + \frac{b_t}{\sigma \sqrt{t}} \, x_t \,.
\end{equation}

\noindent
By separation of the time dimension from the vibrating variable, yields correspondance $f(y_t) = f(t, y_t)$ in bivariate representation, where $y_t$ is a linear application of $x_t$ as per (37). By applying bivariate Taylor expansion to function $f(t, y_t)$ having for support bilinear correspondance between vibrating variable $y_t$ and the source of vibrations $x_t$ (up to second order), leads to:

\begin{multline}
\delta f(t, y_t)  = \frac{\partial f(t, y_t)}{\partial t} \delta t + \frac{\partial f(t,y_t)}{\partial y_t} \delta y_t + \frac{\partial f^2 (t,y_t)}{\partial y_t \partial t} \delta y_t \delta t +  \frac{1}{2} \frac{\partial f^2 (t,y_t)}{\partial t^2} \delta t^2  + \\
+ \frac{1}{2}\frac{\partial f^2(t,y_t)}{\partial y_t^2} \delta y_t^2  +...  \,.
\end{multline}

\noindent
where $\delta f(t,y_t) = f(t+ \delta t, y_t + \delta y_t) - f(t, y_t)$ and $\delta y_t = \frac{b_t}{\sigma  \sqrt{t}} \delta x_t$ by bilinearity between $x_t$ and $y_t$.

\vspace{2mm}
\noindent
Let us rewrite (38) by expanding $\delta y_t$, leading to:

\begin{multline}
\frac{\delta f(t, y_t)}{\delta t}  = \frac{\partial f(t, y_t)}{\partial t} + \frac{b_t}{\sigma \sqrt{t}} \frac{\partial f(t,y_t)}{\partial y_t} \frac{\partial W_t}{ \partial t} + \frac{b_t}{\sqrt{t}}\frac{\partial f^2 (t,y_t)}{\partial y_t \partial t} \delta W_t  +  \frac{1}{2} \frac{\partial f^2 (t,y_t)}{\partial t^2} \delta t  + \\
+ \frac{1}{2} \frac{b_t^2}{t} \frac{\partial f^2(t,y_t)}{\partial y_t^2} \frac{\partial W_t}{ \partial t} \, \delta W_t +...  \,.
\end{multline}

\noindent
where $\frac{\partial W_t}{\partial t}$ is the partial derivative of a Wiener process in integral form.

\vspace{2mm}
\noindent
The terms in $\delta W_t$ and $\delta t$ in (39) vanish and by expressing $\frac{\partial W_t}{\partial t}=\frac{1}{2 t^{1/2}} Z$ we get:

\begin{equation}
\frac{\delta f(t, y_t)}{\delta t}  = \frac{\partial f(t, y_t)}{\partial t} +\frac{1}{2} \frac{b_t}{t} \frac{\partial f(t,y_t)}{\partial y_t} \, Z  \,.
\end{equation}

\vspace{1mm}
\noindent
In context of (36), by some approximation we have $\frac{\partial f(t, y_t)}{\partial y_t} \sim t \, \frac{\partial f(t,y_t)}{\partial t}$ as the integral of ${e}^{y_t}$ over a time interval $t$ is the average value of the function over that interval multiplied by the width of the interval. 

\vspace{2mm}
\noindent
The application of the $\partial /\partial t$ operator to non-homogeneous It{\^o}'s processes as in (36), results in an additional term $\Upsilon = \frac{1}{2} \frac{b_t}{t} \frac{\partial f(t,y_t)}{\partial y_t} \, Z$ as per (40). By crossover of $\Upsilon$ with the other terms  in (41) leads to additional covariances. As a simplification, in the below we applied a single adjustment to $\frac{\partial f(t, y_y)}{\partial t}$, expressed as a factor product $\text{\textgamma}=\sqrt{1 + 1/2 \, \sigma^2 }$. 

\vspace{2mm}
\noindent
As applied to SDE (36) with adjustment factor $\text{\textgamma}=\sqrt{1 + 1/2 \, \sigma^2 }$, leads to: 

\begin{equation}
\frac{\partial y_t}{\partial t}=  \mu +\frac{1}{2}\sigma^2 + \text{sign}(\theta) \, \frac{1}{2}\sigma^2 \, \text{\textgamma}\, {e}^{y_t} + \frac{1}{2} \sigma t^{-1/2} Z \,.
\end{equation}

\vspace{1mm}
\noindent
As from $y_t= m_t + \sqrt{v_t} \, Z$ where $Z \sim \mathcal{N}(0, 1)$, the variance of (41) is as follows:

\begin{equation}
Var \left(\frac{\partial y_t}{\partial t} \right)=  \frac{1}{4} \sigma^4 \text{\textgamma}^2 Var\left( {e}^{m_t + \sqrt{v_t} Z}\right) + \frac{1}{4} \frac{\sigma^2}{t} + \text{sign}(\theta)\frac{1}{2} \frac{\sigma^3}{ t^{1/2}} \text{\textgamma} \, Cov \left( {e}^{m_t+\sqrt{v_t} Z }, Z\right)\,.
\end{equation}

\noindent
The identity $Var(X+Y) = Var(X) + Var(Y) + 2 Cov(X,Y)$ was invoked in (42).

\vspace{1mm}
\noindent
For a lognormal random variable $u={e}^x$ where $x \sim \mathcal{N}(\mu, \sigma^2)$, we have $Var(u)=\left({e}^{\sigma^2} -1 \right) {e}^{2 \mu + \sigma^2}$, leading to:

\begin{equation}
Var\left( {e}^{m_t + \sqrt{v_t} Z}\right) = \left({e}^{v_t} -1 \right) {e}^{2 m_t + v_t}\,.
\end{equation}

\noindent
The valuation of the covariance term in (42), invokes formula $Cov(X,Y)= \mathbb{E} (XY) -  \mathbb{E} (X) \mathbb{E} (Y)$. As $\mathbb{E} (Z) =0$, the covariance term equals  $\mathcal{I} = \mathbb{E} \left( Z {e}^{m_t+\sqrt{v_t} Z}\right)$.

\noindent
We have:

\begin{equation}
\mathcal{I} = \int_{-\infty}^{\infty} x {e}^{m_t+\sqrt{v_t} x}\frac{1}{\sqrt{2 \pi}} {e}^{-\frac{1}{2} x^2} dx \,.
\end{equation}

\noindent
We can write:

\begin{equation}
\mathcal{I} =\frac{1}{\sqrt{2 \pi}} \int_{-\infty}^{\infty} x {e}^{m_t+\sqrt{v_t} x - \frac{1}{2} x^2} dx \,,
\end{equation}

\noindent
leading to:

\begin{equation}
\mathcal{I} = - \frac{1}{\sqrt{2 \pi}} \int_{-\infty}^{\infty}(\sqrt{v_t} -x) {e}^{m_t+\sqrt{v_t} x - \frac{1}{2} x^2} dx + \sqrt{\frac{v_t}{2 \pi}}  \int_{-\infty}^{\infty} {e}^{m_t+\sqrt{v_t} x - \frac{1}{2} x^2} dx  \,.
\end{equation}

\noindent
The first term of integral $\mathcal{I}$ is equals to zero, by asymptotic convergence as both branches tend to infinity, leading to:

\begin{equation}
\mathcal{I} = \sqrt{\frac{v_t}{2 \pi}}  \int_{-\infty}^{\infty} {e}^{m_t+\sqrt{v_t} x - \frac{1}{2} x^2} dx  \,.
\end{equation}

\noindent
As we have $m_t+\sqrt{v_t}x-\frac{1}{2}x^2= -\frac{1}{2} \left( x - \sqrt{v_t}\right)^2 +m_t + \frac{1}{2} v_t$, we can write:

\begin{equation}
\mathcal{I} = \sqrt{v_t} {e}^{(m_t + \frac{1}{2} v_t)}  \int_{-\infty}^{\infty}\frac{1}{\sqrt{2 \pi}} {e}^{-\frac{1}{2}(x-\sqrt{v_t})^2} dx  \,.
\end{equation}

\noindent
The integrant in (48) is the density function of a Gaussian variable with mean $\sqrt{v_t}$ and unit variance. Its integral over the real domain is equal to one. We get:

\begin{equation}
\mathcal{I} = \sqrt{v_t} \, {e}^{(m_t + \frac{1}{2} v_t)} \,.
\end{equation}

\noindent
Hence:

\begin{equation}
Var \left(\frac{\partial y_t}{\partial t} \right) =  \frac{1}{4} \text{\textgamma}^2 \sigma^4 \left({e}^{v_t} -1 \right) {e}^{2 m_t + v_t}+ \frac{1}{4} \frac{\sigma^2}{t} + \text{sign}(\theta) \,\frac{1}{2} \frac{\sigma^3}{ t^{1/2}}\,\text{\textgamma} \,  \sqrt{v_t} \, {e}^{(m_t + \frac{1}{2} v_t)}  \,.
\end{equation}

\noindent
By coupling of (34) with (50), i.e. $y_t$ and $\widetilde{y_t}$ expressing the same SDE, we get:

\begin{equation}
\frac{\partial v_t}{\partial t} = \sqrt{\frac{v_t \sigma^2}{t}+ v_t \sigma^4 \text{\textgamma}^2 \left({e}^{v_t} -1 \right) {e}^{2 m_t + v_t}+ \text{sign}(\theta) \, 2 \frac{\sigma^3}{  t^{1/2}}\, \text{\textgamma}\, v_t^{3/2} \, {e}^{(m_t + \frac{1}{2} {v_t})}}  \,,
\end{equation}

\noindent
with initial conditions $v_0=0$ and $v_t^{\prime}|_{t=0}=\sigma^2$ and $\text{\textgamma}=\sqrt{1 +1/2 \sigma^2}$. Eq. (51) is an approximation of the differential equation describing the variance of the process $y_t$, at any instant of time $t \geq 0$ and where $\theta$ is real.  

\section{The statistical expectation of $f(x)$ with $x \sim \mathcal{N}(\mu, \sigma^2)$ by thin-tile integration method}

Thin-tile integration is an efficient method to compute the statistical expectation of a continuous function $f\colon \mathbb{R} \rightarrow \mathbb{R}$ of a Gaussian random variable $x \sim \mathcal{N}(\mu, \sigma^2)$,  where $\mu$ is the mean and $\sigma$ the standard deviation. This method consists of a non-uniform grid spacing built as a continuum of thin tiles such as shown in fig. 1, which further benefits from the symmetry of the Gaussian distribution.

\begin{figure}[h]
\includegraphics[width=7.7cm, height=4.6cm]{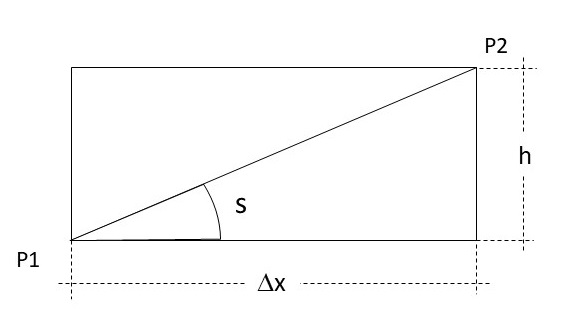}
\centering
\caption{Thin tile connecting two points P1 and P2 of the density function of the Gaussian distribution, where $h$ is the height and $\Delta x$ the width.}
\centering
\end{figure}

\vspace{1mm}
\noindent
The rationale involved by the method is to dynamically determine the width of grid elements from adjacent tiles, connecting points of the density function of the Gaussian distribution by their edges, and where height $h$ is fixed. By considering squarish tiles of sides $h \times h$, the relation betwen the number of pairs of tiles $N$ and $h$ is expressed as $h=\sqrt{1/(2 N)}$. The slope coefficient for a tile is defined as $s = h/\Delta x$, yielding a tile area expressed as $A=h^2 / s$. The tiles are placed two-by-two in a symmetrical fashion on both sides of the Gaussian density function, in ascending order starting from the mode of the curve and prolonged to the tails. As the overall area under the density function between  both extremum spanned by the first $n$ pairs of tiles is described by relation $A_n=1 - 2 \, \Phi(\frac{\mu-x_n}{\sigma})$,  where $\Phi$ is the cumulative density function of the standard Gaussian distribution $\mathcal{N}(0, 1)$ and $x_n$ the coordinate to the right of the surface spanned by the tiles, the slope coefficient impacts the cumulative area $A_n$ by the application of tiny increments $\Delta A_n = 2 \, \frac{h^2}{s_n}$ representing the area of the tiles. By the squared tile rule, say the slope coefficient used to determine the area at the $n^{th}$ pair of tiles disposed on the curve, is set to be equal to the derivative of the density function floored to one, leading to:

\begin{equation}
s_n = \max \left(1.0, \left[ \frac{\mu-x_{n-1}}{\sigma}  \right] \, \varphi(x_{n-1}) \right) \,,
\end{equation}  

\noindent
where $ \varphi(x)=\frac{1}{\sigma\, \sqrt{2\, \pi}} \, \exp\left( - \frac{(x-\mu)^2}{2\, \sigma^2}\right)$ is the density function of the Gaussian distribution  $\mathcal{N}(\mu, \sigma^2)$ of mean $\mu$ and standard deviation $\sigma$, where $x_{n-1}$ represents the rightmost extremum at the $n^{th}-1$ pair of tiles disposed on the curve and where $x_0=\mu$ is the mode. 

\begin{figure}[h]
\includegraphics[width=8.5cm, height=3.5cm]{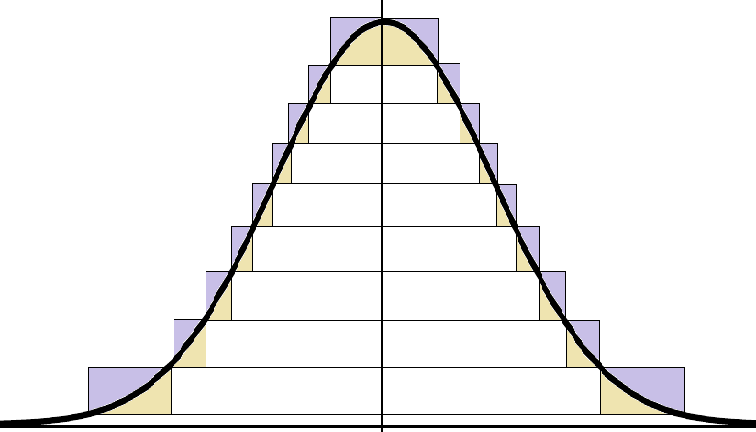}
\centering
\caption{Vault formed by the disposition of thin-tiles on the Gaussian density function defining contours of height $n h = \varphi(\mu)- \varphi(x_n)$ at the $n^{th}$ pair of tiles, where $\varphi(x)$ stands as the Gaussian density function $\sim \, \mathcal{N}(\mu, \sigma^2)$ of mean $\mu$ and standard deviation $\sigma$ positive.}
\centering
\end{figure}

\noindent
From the above, the coordinate referring to the $n^{th}$ pairs of tiles on the right side of the curve, by surface is as follows:

\begin{equation}
x_n = \mu - \sigma \,\, \Phi^{-1}\left(\frac{1-A_n}{2} \right) \,,
\end{equation} 

\noindent
where $\Phi^{-1}$ is the inverse cumulative density function of the standard Gaussian distribution $\mathcal{N}(0, 1)$ and $A_n$ the cumulative area under the density function of say standard tiles of area $h \times h$. Thinner convergence at the tails is achieved with a scheme involving non-standard tiles as per contours as shown in fig. 2.

\vspace{1mm}
\noindent
For each step indexed by $n=1,...,N-1 \in \mathbb{N^{\ast}}$ corresponding to a pair of tiles disposed on Gaussian distribution, yields an observable of the function of the Gaussian variable, which by arithmetic mean from points disposed symmetrically around the mode, leads to:

\begin{equation}
f_n = \frac{1}{4} \, \left[ f(x_n) + f(x_{n-1}) + f(\mu-x_n) +f(\mu-x_{n-1}) \right] \,,
\end{equation} 

\noindent
as an observable of $f(x)$ weighted by the incremental area $\Delta A_n$ corresponding to the $n^{th}$ pair of tiles disposed on the curve in the standard scheme. We finally compute the weighted averge of the $N-1$ observations of $f(x)$ representing the statistical expectation of $f(x)$ under natural probabilities, where $x \sim \mathcal{N}(\mu, \sigma^2)$ as a Gaussian distribution of mean $\mu$ and standard deviation $\sigma$.

\section{Numerical results}

The MGF of the lognormal distribution with parameter $\mu$, $\sigma$ and where $\theta$ is real, is expressed as $M(\theta)=\mathbb{E}(f_1)$ where the shifted process $\widetilde{f_t}$ is coming from SDE (22, 23). The parametric functions $m_1$ and $v_1$ are evaluated by integration of differential equations (31) and (51) over a unit time interval $\Delta T=[0,1]$. 

\vspace{2mm}

\noindent
We solve these integrals by discretisation over the domain of integration $\Delta T$, introducing small time steps $\delta t$. We start the calculation from time $t_0=0$ and iteratively compute $\frac{\partial m_t }{\partial t}$ and $\frac{\partial v_t}{\partial t}$ using piecewise linear segments, leading to the valuation of $m_t$ and $v_t$ at the next time step until reaching $t_1=1$. The numerical scheme consists of:

\begin{equation}
m_{i+1} = m_{i} + \frac{\partial m_t}{\partial t}\Big|_{i} \delta t \,,
\end{equation}

\noindent
and

\begin{equation}
v_{i+1} = v_{i} + \frac{\partial v_t}{\partial t}\Big|_{i} \delta t \,,
\end{equation}

\noindent
at each iteration.

\vspace{2mm}

\noindent
The initial conditions are given by $m_0=\ln (|\theta|)$, $v_0=0$ and $v_t^{\prime}|_{t=0}=\sigma^2$, by convergence of the variance of the process towards $t \, \sigma^2$ as $t$ tends towards zero. Once we get the endpoint $m_1$, the estimator of the MGF of the lognormal is given by $\hat{M}(\theta)=e^{\text{sign}(\theta)\, e^{m_1}}$. This is the approach for the valuation of the MGF of the lognormal by the stochastic approach based on zero-entropy principle, i.e. $f_t$ shifted to $\widetilde{f_t}$ under non-vibrating variable, and reciprocal by symmetry.


\noindent
As a reference, the MGF of the lognormal distribution as given by the Laplace transform of the lognormal in \cite{Asmussen2016}, is expressed as follows:

\begin{equation}
\hat{M}_{\mathcal{L}}(\theta) \approx \frac{\exp \left( - \frac{W^2(-\theta \sigma^2 {e}^{\mu})+ 2 W(-\theta \sigma^2 {e}^{\mu})}{2 \sigma^2} \right)}{\sqrt{1+ W(-\theta \sigma^2 {e}^{\mu})}} \,,
\end{equation}

\noindent
where W is the Lambert-W function defined as the inverse of $f(w) = w \, \exp(w)$.

\vspace{2mm}
\noindent
The valuation of the MGF of the lognormal are summarized in tables 1, 2 and 3 for a range of $\theta$ values. Parameter $\mu$ was set to zero in all three tables, whereas standard deviation set to $\sigma=0.1$ for table 1, $\sigma=0.0625$ for table 2 and $\sigma=1.0$ for table 3. With regard to the zero-entropy stochastic approach, 2'000 equidistant time steps were used part of the discretisation algorithm. The accuracy was set to $1.0\e{-6}$ for the Lambert function as part of the Asmussen, Jensen and Rojas-Nandayapa method. For thin-tile integration, standard tiles with height set to $h=0.0025$ were used, i.e. $N=80,000$. In contrast, plain vanilla Monte Carlo simulation for the estimation of the MGF lognormal required a sample size of about $100$ millions of observations to achieve commensurate accuracy level. Note that in table 1 where $\theta$ is positive, the variance was initialised to $v_0=\sigma^2$ when using the approximate factor $\text{\textgamma}$, which is a side effect for not using proper covariances of $\Upsilon$ with the other terms in (41). As a hint, the approximation $\Upsilon \approx \frac{1}{\sqrt{2}} \, \sigma \, \frac{\partial f(t, y_t)}{\partial t}$ yields an additional term $\frac{1}{2^{3/2}} \, \sigma^5 \, \sqrt{v_t} \, {e}^{2\, m_t + v_t}$ from its covariance with non-homogeneous term $\frac{1}{2} \sigma^2 \, {e}^{y_t}$ in (41). Whenever $\theta$ is negative as in table 2 and 3, this variable is properly initialised to $v_0=0$.

\begin{table}[!ht]
\centering
\caption{Table for MGF of the lognormal when $\theta$ is positive with $\sigma = 0.1$}
\begin{tabular}{| p{4.0cm} | c |c|c|c|c|c|}
  \hline 
$\theta$ & 0.1 & 0.3 & 0.5 & 1.0  & 1.2 \\
  \hline \hline
 \hline
Monte Carlo simulation & 1.105779 & 1.352510 & 1.654955 & 2.745936  & 3.365014 \\
\hline
Stochastic approach based on zero-entropy principle & 1.105780 & 1.352506 & 1.654957 & 2.745994  & 3.365088 \\
\hline
Thin-tile integration method & 1.105781 & 1.352509 & 1.654966 & 2.745978 & 3.364940 \\
\hline
Asmussen, Jensen and Rojas-Nandayapa approximation & 1.105780 &1.352504 & 1.654957 & 2.745950  & 3.364990\\
\hline
\end{tabular}
\vspace{6mm}

\centering
\caption{Table for MGF of the lognormal when $\theta$ is negative with $\sigma=0.0625$}
\begin{tabular}{| p{3.8cm} | c |c|c|c|c|c|}
  \hline 
$\theta$ & -0.5 & -1.0 & -2.0 & -4.0 & -8.0   \\
  \hline \hline
 \hline
Monte Carlo simulation  & 0.606235 & 0.367884 & 0.135863 & 0.018744  &   0.000373\\
\hline
Stochastic approach based on zero-entropy principle  & 0.606234 & 0.367879 & 0.135863 & 0.018746 & 0.000373 \\
\hline
Thin-tile integration method & 0.606235 & 0.367880 & 0.135862 & 0.018744 & 0.000373   \\
\hline
Asmussen, Jensen and Rojas-Nandayapa approximation  & 0.606235 & 0.367880 & 0.135862 & 0.018744 & 0.000373   \\
\hline
\end{tabular}
\end{table}

\begin{table}[ht]
\centering
\caption{Table for MGF of the lognormal when $\theta$ is negative with $\sigma=1.0$}
\begin{tabular}{| p{3.8cm} | c |c|c|c|c|c|}
  \hline 
$\theta$ & -0.5 & -1.0 & -2.0 & -4.0 & -8.0   \\
  \hline \hline
 \hline
Monte Carlo simulation  & 0.561707 & 0.381729 & 0.216326  & 0.098069 & 0.034274  \\
\hline
Stochastic approach based on zero-entropy principle  & 0.560233 & 0.367879 & 0.238030 &  0.159668  & 0.118724 \\
\hline
Thin-tile integration method  & 0.561708 & 0.381755 & 0.216305 & 0.098046 & 0.034264   \\
\hline
Asmussen, Jensen and Rojas-Nandayapa approximation  & 0.561717 & 0.381752 & 0.216304 & 0.098042 & 0.034267   \\
\hline
\end{tabular}
\end{table}

\section{Conclusion}

Thin-tile integration and the Laplace transform of the logonormal by Asmussen, Jensen and Rojas-Nandayapa are in good agreement, giving values for the lognormal MGF in all three settings: positive $\theta$ values,  negative $\theta$ values, and high volatilities, i.e. $\sigma=1.0$, at an accuracy of $5$ to $6$ digits after the decimal point, providing a valuable benchmark for the stochastic approach discussed below.

The stochastic approach by the application of It{\^o}'s calculus to the lognormal MGF, is based on the so-called zero-entropy principle and the symmetry of an action and reciprocal. While the stochastic approach yields lognormal MGF values in pretty good agreement with the aforementioned methods for negative $\theta$ values, matching the benchmark with an accuracy of about $6$ digits after the decimal (see table 2), tiny differences in the order of 0.1 to 1.0 basis points were obtained in table 1 for positive $\theta$ values, and larger deviations in table 3 when applied to vibrations carrying higher volatilities, i.e. $\sigma$ set to one. Although variations observed in table 1 and 2 could be attributed to numerical imprecisions, the slight departures from the benchmark occuring at higher volatilities as seen in table 3 are of statistical significance. As an explanation for these small differences, present study provides support for the partial offset by the reciprocal of an action as applied to SDEs carrying vibrations, resulting from an asymmetry under zero-entropy principle or non-linearity of the $\mathbb{E}_{\ast}$ operator of probability space $\mathcal{E}_{\cdot}$ corresponding to statistical expectations in epsilon probability measure. Notwithstanding, the foregoing does not preclude from other sources of inaccuracies such as crossover covariances from $\Upsilon$ with other terms in (41). 

Moreover, non-Gaussianity of the underlying process when extended to the complex domain as seen in §3.2, is preventing the lognormal MGF from the stochastic approach to be applicable to complex numbers, a prerequisite for the characteristic function. As puzzling as this may be, non-Gaussianity of the underlying process when extended to complex numbers is an aspect having a connection with application of Euler's formula and the bivariate structure of the underlying stochastic differential equations. Finally, the asymmetry from the partial offset by the reciprocal of an action as applied to systems carrying vibrations, is a principle having potential applications in other fields such as quantum theory, e.g. asymmetries of harmonic oscillators, zero-point vibration preventing liquid helium from freezing at atmospheric pressure, etc.

\newpage


\bibliographystyle{plain}
\bibliography{lognormal}   

%
%

\end{document}